\documentclass[12pt]{article}

\usepackage[a4paper]{geometry}
\usepackage{ajc-FR}
\usepackage{amsmath, amsfonts, latexsym}
\usepackage{graphicx}
\usepackage{float}
\usepackage{hyperref}
\usepackage{placeins}

\Volume{XY(Z)}
\Year{2020}
\firstpage{1}
\lastpage{5}
\Revised{24 November 2020}
\runninghead{\copyright \, F. Rowley 2020\,/\,NLBWSP DRAFT arXiv version, pages 1 -- 5}
\numberwithin{equation}{section}

\newtheorem{theorem}{Theorem}[section]

\newenvironment*{proof}
{\begin{list}{}{\setlength{\leftmargin}{0em}\setlength{\rightmargin}{0em}}
\item[] {\sc Proof:}} {\hfill$\Box$
\end{list}}

\begin{document}

\title{\bf New lower bounds for\\ weak Schur partitions}
\author{Fred Rowley\thanks{formerly of Lincoln College, Oxford, UK.}}
\Addr{West Pennant Hills, \\NSW, Australia.\\{\tt fred.rowley@ozemail.com.au}}

% \date{\dateline{submission date}{acceptance date}\\
% \small Mathematics Subject Classifications: comma separated list of
% MSC codes available from http://www.ams.org/mathscinet/freeTools.html}

\date{\dateline{24 November 2020}{DD Mmm CCYY}}

\maketitle

% Papers must include an abstract. The abstract should consist of a
% succinct statement of background followed by a listing of the
% principal new results that are to be found in the paper. The abstract
% should be informative, clear, and as complete as possible. Phrases
% like "we investigate..." or "we study..." should be kept to a minimum
% in favour of "we prove that..."  or "we show that...".  Do not
% include equation numbers, unexpanded citations (such as "[23]"), or
% any other references to things in the paper that are not defined in
% the abstract. The abstract will be distributed without the rest of the
% paper so it must be entirely self-contained.

\begin{abstract}

This paper records some apparently new results for the partition of integer intervals $[1, n]$ into weakly sum-free subsets.  These were produced using a method closely related to that used by Schur in 1917. 

New lower bounds can be produced in this way for partitions of unlimited size.  The asymptotic growth rate of the lower bounds, as the number of subsets increases, cannot be less than the same growth rate for strongly sum-free partitions, and so exceeds 3.27.  

Specific results for partitions into a 'small' number of subsets include $WS(6) \ge 642$, $WS(7) \ge 2146$, $WS(8) \ge 6976$, $WS(9) \ge 21848$, and $WS(10) \ge 70778$.

% keywords are optional
\bigskip\noindent
%\bigskip\noindent \textbf{Keywords:} weak Schur partition, sum-free sets.
\end{abstract}
\bigskip
\small DRAFT \copyright Fred Rowley  November 2020.
\bigskip

\section{Introduction}
It is natural to ask how the most successful methods used to colour linear triangle-free graphs (or equivalently, to construct strong Schur partitions) might be modified so as to permit the construction of larger \textbf{weak} Schur partitions.

The author's previous experience indicates that the most successful constructions for 'small' triangle-free graphs can be characterised as special cases having particular unique attributes.  Smaller graphs have typically been derived by a range of exhaustive or partial search strategies: and have then been combined or extended by methods involving various forms of 'compounding'.  Compound graphs may be derived using periodically repetitive structures (translations) and/or reflections.  These techniques usually succeed by vastly reducing the size of the difference sets derived from the subsets comprised in the colouring.  

Any strong Schur partition is also a weak Schur partition, so both the size of any maximal weak partitions into $r$ subsets, and their ultimate growth rate as $r$ increases, cannot  be less than in the strong case. Previous papers, including \cite{FR-GLRGC}, have demonstrated that the ultimate growth rate for strong Schur partitions, as the number of colours $r$ increases, exceeds 3.27.  This author has seen evidence that some constructions become much more difficult when the ratio approaches $(3+\sqrt{13})/2$, which is a little over 3.3.

One immediate observation, when attempting to construct weak Schur partitions, is that translations or reflections are much less useful.  If a pair $(x, 2x)$ exists in a single subset $S_i$ in a weak Schur partition, it is clearly not possible in general to feature either of the pairs $(x+a, 2x+a)$ or $(m-x, m-2x)$ in the same subset, since in each case the difference is $x$, and $x \in S_i$.  

Some recent work of this topic has succeeded in increasing the known lower bounds by sidestepping these constraints using various algorithms and search constraints - see, for instance, \cite{Rafili}. So far, however, it has apparently not proved possible to demonstrate in this way an infinite sequence of weak Schur partitions with a growth rate above 3, which does not consist simply of strong partitions.  

This paper now provides two such sequences.  Although one might rightly say that the partitions in these sequences are 'almost' strong partitions, that may simply indicate that there is room for more imaginative constructions.  

Notation is defined in section 2.  

In section 3, it is proved that, starting from a single specific graph, a series of graphs can be constructed, giving improved values for $WS(s)$ applicable for all $s \ge 6$. Numerical lower bounds are shown for $1 \le r \le 10$.  

In section 4, some very brief conclusions are drawn. 

\section{Definitions and Notation}

In this paper: 

The set of integers $\{1, 2, 3, \dots, n\}$ is sometimes written as a \emph{closed integer interval} $[1, n]$. A partition $p$ of $[1, n]$ into $r$ subsets may be denoted by $p(r; n)$.

If an interval $S = [1, n]$ can be partitioned into $r$ disjoint non-empty subsets $S_i$ for $i = 1, 2, 3, \dots, r$, where no subset contains three \textbf{distinct} integers $a, b, c$, such that $a+b=c$, then each such subset is \emph{weakly sum-free} and that partition is a \emph{weak Schur partition}. The order of the set $S$ is clearly $n$, and is also referred to as the \emph{order of the partition}. 

A pair of positive integers $(a, 2a)$ is referred to as a \emph{weak pair}.

For any $r$, $WS(r)$ is the maximum value of $n$ such that a weak Schur partition $p(r; n)$ exists.  $WS(r)$ is known as the \emph{weak Schur number}, and its existence is established by Ramsey's Theorem. 

\section{Construction of Weak Schur Partitions}

\begin{theorem}
  \label{Thm:C-Thm1}
(Construction Theorem)\\
If there is a strong Schur partition of the integers $[1, m]$ into $r$ subsets, then there is a weak Schur partition of $[1, 4m+2]$ into $r+1$ subsets; and a weak Schur partition of $[1, 13m+8]$ into $r+2$ subsets.    
\end{theorem}

The theorem depends on two very simple constructions, which are closely related to that used by Schur in $\cite{Sch1}$.   

\begin{proof}

As stated above, the repetition or reflection of 'weak pairs' $(a, 2a)$ within a prototype partition, into an extended partition, is not useful in the general case.  The first construction minimises this problem by relying on a sequence of partitions, each of which has only the single 'weak pair', $(1, 2)$ in one of its subsets.  No other weak pairs are involved.  

The first construction takes as its first 'prototype' the following partition of order 6:

$S_1 = \{1, 2, 6\}$\\
$S_2 = \{3, 4, 5\}$.

We also assume the existence of a strong Schur partition $q(r; m)$.

First, $S_1$ is extended to include an arithmetic series with difference 4, and so becomes $T_{r+1} = \{1, 2, 6, 10, \dots, 4m+2\}$.

$S_2$ is used to construct $m$ 'translates' of the numbers 3, 4 and 5 -- giving us $m$ distinct subsets of $[1, 4m+2]$, each of the form $T_i = \{4i-1, 4i, 4i+1\}$.  It can be seen that within any one of these subsets $T_i$, the absolute differences are either 1 or 2 and so are members of $S_1$.  Therefore each subset $T_i$ is sum-free.  

We then form the remaining subsets of the new partition, $S_1, S_2, \dots , S_r$ by taking the unions of all the subsets $T_i$ whose indices are in the same subset in the strong partition $q(r; m)$.

Let us assume that two distinct subsets $T_i, T_j$ (with $i < j$) are included in the same subset of the new weak Schur partition. If so, any difference between a member of $T_j$ and a member of $T_i$ must be in the interval $[4(j-i)-2, 4(j-i)+2]$.  Any number in this range is always either (a) a member of the subset $T_{j-i}$; or (b) a member of $S_1$.  

Case (a) is the only case we need concern ourselves with.  In that case, the partition $q(r; m)$ would not be sum-free if the subset containing $(j-i)$ were the same as that containing both $i$ and $j$.  Therefore, every subset in the new partition, which was formed by taking the union of the $T_i$, is also strongly sum-free.  Clearly, the subsets in the new partition that consist of unions of the $T_i$ must be of the same cardinality as the subsets in $q(r; m)$ -- i.e. $r$. Therefore, including $T_{r+1}$, we have $r+1$ subsets in the new partition

Lastly we observe that $T_{r+1}$ is weakly sum-free, that the order of the new partition is $4m+2$, and that it is a complete partition of $[1, 4m+2]$. 

We omit many details of the proof of the second construction used in the theorem, which partitions the set $[1, 13m+8]$.  It follows exactly similar lines to the above, but starts from the following partition of $[1, 21]$ into 3 subsets: 

$S_1 = \{1, 2, 4, 8, 21\}$\\
$S_2 = \{3, 5, 6, 7, 18, 19, 20\}$\\
$S_3 = [9, 17]$.

The sets $T_i$ are derived in this case by translation of $[9,17]$. Set $T_{r+1}$ is derived by extending $S_1$ and set $T_{r+2}$ is derived by extending $S_2$, in each case with a period of 13.  We note that only values above 4 are included in the basis for these extensions; and that each partition produced in this way contains exactly four weak pairs.

\end{proof}

It is now simple to deduce that if there is an infinite sequence of strong partitions with an ultimate growth rate of (say) $\gamma$, then there is a corresponding sequence of weak partitions with the same ultimate growth rate.

The orders now available for some 'small' partitions are shown in Table 1 below.  The history and derivation of the smaller weak Schur partitions is well covered in \cite{EMRS}, and orders of weak partitions shown for $1 \le r \le 5$ are from that source. All are believed to be the largest currently available, and the first four have been shown to be maximal.

The orders of weak partitions for $6 \le r \le 10$ are produced by the construction above and believed to exceed the highest values previously published. The orders shown for the strong partitions (on which the constructions are based) are derived from \cite{FR-GLRGC}. \\

\FloatBarrier

\textbf{Table 1 - Orders of largest available weak and strong Schur partitions}

%\begin{table}[!ht]
%  \begin{center}
    % use \includegraphics to import figures 
    \includegraphics{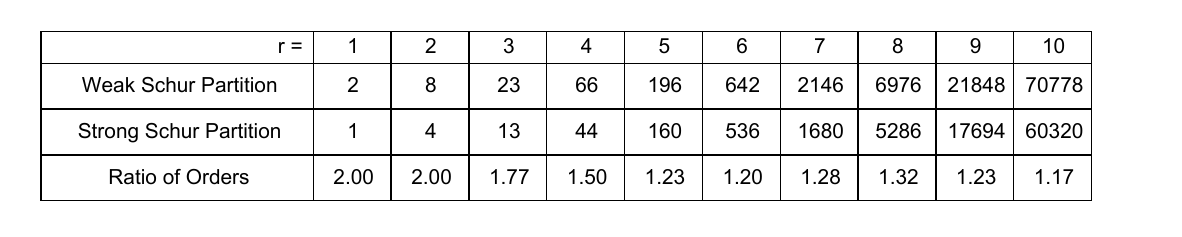}
%  \end{center}
%  \caption{\label{Table:01} Lower Bounds on WS(r)}
%\end{table}

\FloatBarrier

\section{Conclusions}
	
The construction described here is quite simple and effective, but fairly limited. 

All of the partitions demonstrated above have very few weak pairs and so may be said to be only 'trivially weak', with order equal to a fixed multiple of a known strong partition.  As a result, we have not shown that the limiting growth rate in the weak case exceeds that in the strong case.  Although there is room for a lot more work, this author believes that the limiting growth rates may well be finite and equal, and further are quite likely to be bounded by a number well below 4.

Nor does this paper provide a sequence in which every partition certainly exceeds the maximum possible strong partition: although it might do so if the ultimate growth rate in the strong case can later be shown to be less than 4.  For the moment, though, many of the 'small' partitions represent significant improvements over previously demonstrated lower bounds on $WS(r)$.

Despite its limitations, the construction demonstrated in this paper sets a new baseline for constructing infinite sequences of weak Schur partitions in a way that consistently exceeds what is possible in the strong case.

An example of a weak Schur partition $q(6; 642)$ is attached to this pdf as an ancillary file.

% \FloatBarrier
% Use \FloatBarrier if using package "placeins"' - this commmand fixes position of Tables - they cannot float below here.

\medskip
%%%%%%%%%%%%%%%%%%%%%%%%%%%%%%%%%%%%%%%%%%%%%%%%%%%%%%%
 \subsection*{Dedication}
I dedicate this paper to the memory of my very good friend, the late Paul~A.~Stanway, former Exhibitioner of St.~John's College Cambridge.

%%%%%%%%%%%%%%%%%%%%%%%%%%%%%%%%%%%%%%%%%%%%%%%%%%%%%%%
% \bibliographystyle{plain} 
% \bibliography{myBibFile} 
% If you use BibTeX to create a bibliography
% then copy and past the contents of your .bbl file into your .tex file

\end{document}